\documentclass[article]{elsarticle}

\usepackage{color}
\usepackage{bm}
\usepackage{multirow}
\usepackage{lineno,hyperref}
\modulolinenumbers[5]

\usepackage{amsopn}

\usepackage{lipsum}
\usepackage{amsfonts}
\usepackage{amsmath}
\usepackage{amsthm}
\usepackage{graphicx}
\usepackage{epstopdf}

\usepackage{algorithmic}
\usepackage{algorithm}
\ifpdf
  \DeclareGraphicsExtensions{.eps,.pdf,.png,.jpg}
\else
  \DeclareGraphicsExtensions{.eps}
\fi

\begin{document}

\begin{frontmatter}

\title{A two level method for isogeometric discretizations}

\author{\'Alvaro P\'e de la Riva \fnref{alvaro}}
\fntext[alvaro]{IUMA and Applied Mathematics Department, University of Zaragoza, Spain,(apedelariva@unizar.es).}
\author{Carmen Rodrigo \fnref{carmen}} 
\fntext[carmen]{IUMA and Applied Mathematics Department, University of Zaragoza, Spain,(carmenr@unizar.es).}
\author{Francisco J. Gaspar \fnref{francisco}} 
\fntext[francisco]{IUMA and Applied Mathematics Department, University of Zaragoza, Spain,(fjgaspar@unizar.es).}
%


\begin{abstract}
Isogeometric Analysis (IGA) is a computational technique for the numerical approximation of partial differential equations (PDEs). This technique is based on the use of spline-type basis functions, that are able to hold a global smoothness and allow to exactly capture a wide set of common geometries. The current rise of this approach has encouraged the search of fast solvers for isogeometric discretizations and nowadays this topic is full of interest. In this framework, a desired property of the solvers is the robustness with respect to both the polinomial degree $p$ and the mesh size $h$. For this task, in this paper we propose a two-level method such that a discretization of order $p$ is considered in the first level whereas the second level consists of a linear or quadratic discretization. On the first level, we suggest to apply one single iteration of a multiplicative Schwarz method. The choice of the block-size of such an iteration depends on the spline degree $p$, and is supported by a local Fourier analysis (LFA). At the second level one is free to apply any given strategy to solve the problem exactly. However, it is also possible to get an approximation of the solution at this level by using an $h-$multigrid method. The resulting solver is efficient and robust with respect to the spline degree $p$. Finally, some numerical experiments are given in order to demonstrate the good performance of the proposed solver. 
 
 
\end{abstract}


\end{frontmatter}


\section{Introduction}

\label{sec:intro}

The IGA technique was firstly introduced by Tom Hughes et al. in the seminal paper \cite{Hughes_CMAME2005} in order to integrate the finite element method (FEM) with the computer aided geometric design. This analysis consists of using spline-type basis functions for the representation of the physical domain, as well as for the numerical approximation of the solution of the PDEs. These functions are globally smooth providing up to $\mathcal{C}^{p-1}$ continuity of the solution, where $p$ denotes the polynomial degree.

Given that the isogeometric discretizations of PDEs yield stiffness matrices whose bandwith becomes wider as $p$ is increased, the search of a robust solver with respect to the spline degree $p$ is not an easy task. However, it is of great interest to obtain efficiently the solution of isogeometric discretizations when high spline degrees are considered. Firstly, in~\cite{Collier2012,Collier2013} a study of the computational efficiency of direct and iterative solvers for IGA, respectively, was performed, and since then, the design of iterative solvers for isogeometric discretizations has attracted a lot of attention.
For example, a multilevel BPX-preconditioner was developed in~\cite{Buffa2013} for isogeometric analysis. Beir\~ao da Veiga et al. analyzed overlapping Schwarz methods for IGA in \cite{Beirao2012}, whereas in \cite{Beirao2013} they studied BDDC preconditioners by introducing appropriate discrete norms. 
Algebraic multilevel iteration (AMLI) methods were applied for the isogeometric discretization of scalar second order elliptic problems in \cite{Gahalaut2013}, and preconditioners based on fast solvers for the Sylvester equation were proposed in \cite{Sangalli2016}.  
In the framework of multigrid techniques, different types of smoothers have been proposed to avoid the troubles encountered by standard relaxation procedures. In \cite{Donatelli2015} a preconditioned Krylov smoother at the finest level was considered and in \cite{Hofreiher2015} the authors proposed a multigrid solver based on a mass matrix smoother. In both cases, an increase in the number of smoothing iterations was needed in order to obtain robustness with respect to the spline degree. To avoid the lack of robustness of the mass smoother, in \cite{Hofreither2017} a new version of such a relaxation including a boundary correction was presented. However, the extension of that version to three dimensions was not clear, and therefore, in \cite{Hofreither2017_b}, the authors proposed a multigrid smoother based on an additive subspace correction technique. In such approach, a different smoother is applied to each of the subspaces: in the regular interior subspace a mass smoother is considered, whereas in the other subspaces they proposed to use relaxations which exploit the particular structure of the subspaces. 
Also $p-$multigrid methods have been applied for solving IGA. In~\cite{Tielen} the authors apply a $p-$multigrid method based on an ILUT (Incomplete  LU  factorization  based  on  a  dual  Threshold strategy) smoother and compare this approach with $h-$multigrid methods based on the same smoother. 
Recently, we have proposed  in~\cite{nuestropaper2018} a very simple robust and efficient geometric multigrid algorithm based on a $V(1,0)-$cycle with overlapping multiplicative Schwarz-type methods as smoothers for solving IGA. The key for the robustness of the algorithm with respect to the spline degree is the choice of larger blocks within the Schwarz smoother when the spline degree grows up.

The main contribution of this work is to propose a robust two-level method for solving a target isogeometric discretizacion of order $p$, such that a linear/quadratic discretization is considered at the second level depending on the parametrization of the physical domain. 
At the first level, we apply only one iteration of a suitable overlapping multiplicative Schwarz method. Then, a restriction operator is constructed via projection of the B-spline basis functions between the corresponding approximation spline spaces of the target degree $p$ and $p=1$ or $p=2$. At this point, the prolongation operator is defined as the adjoint of the restriction operator. For solving exactly the system arising on the second level there exist well-known solution techniques. However, one can also obtain an approximation of the solution at the second level by using few steps of an iterative method. In this work, we propose to apply an $h-$multigrid on the coarse level. More concretely, one single iteration of a $V(1,1)-$cycle that uses a red-black Gauss-Seidel smoother provides very good results. Moreover, a further improvement of the algorithm can be achieved by using a more aggressive coarsening strategy. In addition to reduce the spline degree from $p$ to $1$ or $2$, we propose to coarsen the grid-size $h$ to $2h$ from the first to the second level. 

The proposed two-level method is theoretically studied by a local Fourier analysis. 
This analysis, introduced by Achi Brandt in \cite{Bra77, Bra94}, is the main quantitative analysis for the convergence of multilevel algorithms, and results in a very useful tool for the design of this type of methods. Moreover, in~\cite{rigorousLFA} it has been recently proved that under standard assumptions LFA is a rigorous analysis, providing the exact asymptotic convergence factors of the method. 
LFA has been successfully applied to isogeometric discretizations in \cite{nuestropaper2018} in order to analyze the convergence of an $h-$multigrid method based on multiplicative Schwarz smoothers. In particular, an analysis for any spline degree $p$ and an arbitrary size of the blocks in the smoother is provided in such work. Here, such an analysis is used to choose for each spline degree $p$ the block-size in the multiplicative Schwarz iteration on the first level that provides a robust two-level algorithm. 
Thus, this analysis theoretically supports the convergence of the proposed two-level method. Furthermore, LFA can be also performed to analyze the version of the algorithm in which we approximate the solution at the second level by using an $h-$multigrid method. In that case, a three-grid local Fourier analysis has to be considered in order to take into account the approximation on the second level instead of an exact solve. Finally, again a two-grid LFA is applied to support the enhancement of the algorithm by considering an standard coarsening strategy between the first and the second levels. 

It is not the first time that a two-level method is proposed for high-order discretizations. In the framework of discontinuous Galerkin (DG) methods, in~\cite{twolevel_Ludmil} it was theoretically proved that a suitable additive Schwarz method provides uniform convergence with respect to all the discretization parameters, i.e. the mesh size, the polynomial order and the penalization coefficient appearing in the DG bilinear form. However, in such a work, the block-size of the appropriate additive Schwarz iteration is not provided and here we support its choice by a suitable local Fourier analysis. 

The rest of the paper is structured as follows: In Section~\ref{sec:preliminaries} a brief introduction to the isogeometric analysis is given. Also, we state here a model problem and the basics of B-splines and NURBS. Section~\ref{sec:two-level} is devoted to the presentation of the proposed two-level method. The algorithm, together with its components, are introduced in Section~\ref{sec:two-level1}; the approach in which an $h-$multigrid is applied on the coarse level is explained in Section~\ref{sec:approx}; and finally an improvement of the two-level method based on an aggressive coarsening is presented in Section~\ref{sec:improvement}. In Section~\ref{sec:lfa}, we develop the corresponding LFA in order to support the design of our solver. We perform the LFA for the three versions of the method and we present the corresponding results. In Section~\ref{sec:numerical} two numerical experiments show the good performance of the proposed two-level method. Finally, Section~\ref{sec:conclusions} summarizes the main results of this work and draws some conclusions.

\section{Isogeometric Analysis}\label{sec:preliminaries}

Let us consider the Poisson equation in a $d$-dimensional domain $\Omega=(0,1)^d$ with homogeneous Dirichlet boundary conditions:
\begin{equation}\label{modelproblem}
\left\lbrace\begin{array}{ccc}
-\Delta u  &=& f, \ \  \mbox{in} \quad \Omega, \\ 
u &=& 0, \ \  \mbox{on} \quad \partial \Omega.
\end{array}\right.
\end{equation}
The variational formulation of our model problem \eqref{modelproblem} reads as follows: Find $u \in H_0^1(\Omega)$ such that $a(u,v) = (f,v), \forall v \in H_0^1(\Omega)$, where
$$
a(u,v) = \int_{\Omega} \nabla u \cdot \nabla v \, {\rm d} x, \quad\mbox{and}\quad 
(f,v) = \int_\Omega f v \, {\rm d} x.
$$
The Galerkin approximation of the variational problem is given by: Find $u_h \in V_h$ such that 
\begin{equation} \label{approxvariational}
a(u_h,v_h) = (f,v_h), \quad \forall v_h \in V_h,
\end{equation}
\noindent where $V_h$ is a finite dimensional space. In the isogeometric framework, $V_h$ is a given space of splines whose global smoothness might vary depending on the refinement strategy~\cite{Hughes_CMAME2005}. In this work, we will consider spline spaces of degree $p$ holding maximum continuity, that is, $\mathcal{C}^{p-1}$ regularity. Thus, we will denote the finite dimensional space as $V_h^p$, and the numerical solution as $u_h^p$ to emphasize the dependence on the spline degree $p$. Once we have fixed a finite basis $\{\varphi^p_1, \ldots, \varphi^p_{n_h^p}\}$, ${\rm dim } V^p_h = n_h^p $, the solution of \eqref{approxvariational} can be expanded as a linear combination of the considered spline basis functions. That is, 
$$u_h^p = \sum_{i=1}^{n_h^p} u^p_i \varphi^p_i, $$
\noindent where the coefficients $u^p_i$ are the so-called control points. These coefficients can be computed by solving the linear system $A_p u_p =b_p$, where the stiffness matrix is given by $A_p=(a_{i,j})=(a(\varphi^p_j,\varphi^p_i))_{i,j=1}^{n_h^p}$ and the right hand side is $b_p  = (f,\varphi^p_i)_{i=1}^{n^p_h}$.

\subsection{B-splines}
Firstly, we introduce the univariate B-splines. Hence, we consider the unidimensional counterpart of problem~\eqref{modelproblem}, that is the two-boundary value problem,  
$$
-u''(x) = f(x), \quad x \in \Omega = (0,1), \quad  u(0) = u(1) = 0.
$$ 
Given that B-splines are constructed parametrically, a non-decreasing sequence of real numbers called knot vector is required to state the parameter space. Since we are interested in a spline space of degree $p$ with maximum smoothness, we consider a uniform partition of the interval $(0,1)$ into $m$ subintervals $I_i = ((i-1)h,ih)$, $i = 1, \ldots, m$, with $h=1/m, \; m\in{\mathbb N}$. Thus, let $\Xi_{p,h}$ be an uniform and open knot vector as follows
$$\Xi_{p,h} = \{\xi_1 = \ldots = \xi_{p+1} = 0 < \xi_{p+2} < \ldots \xi_{p+m} < 1 = \xi_{p+m+1} = \ldots
= \xi_{2p+m+1}\},$$
\noindent where $\xi_{p+i+1} = i/m, \; i = 0, \ldots, m$. Given such an open knot vector, we define the spline space of degree $p \geq 1$ with maximum continuity as follows:
\begin{equation}\label{splinespace}
\mathcal{S}_{h}^p(0,1) = \{ u^p_h \in C^{p-1}(0,1) : u_h^p |_{I_i} \in {\mathbb P}^p, i = 1, \ldots, m, u_h^p(0) = u_h^p(1) = 0\}, 
\end{equation}
\noindent where $C^{p-1}(0,1)$ is the space of all $p-1$ times continuously differentiable functions on $(0,1)$, and ${\mathbb P}^p$ is the space of all polynomials of degree less than or equal to $p$. 
The dimension of the space $\mathcal{S}_{h}^p(0,1)$ is $p+m-2$, and the set of basis functions $\left\{ N_i^p \right\}_{i=2}^{p+m-1}$ of this space is defined recursively by the Cox-de-Boor formula (see~\cite{deboor}), starting with $p=0$ (piecewise constants). 
For the case $p=0$, the constant splines are step functions with support on the corresponding knot span. That is, for $i=1,\ldots,m+2p$,
\begin{equation}
\label{constantspline}
N_{i}^0(\xi)=\left\lbrace \begin{array}{ll}
1 & \textrm{if } \xi_{i} \leq \xi < \xi_{i+1}, \\
0 & \textrm{otherwise.} \\
\end{array}\right.
\end{equation}	
Then, for every pair $(k,i)$ such that $1\leq k \leq p$, $1 \leq i \leq m+2p-k$, the basis functions $N_{i}^k : [0,1] \rightarrow {\mathbb R}$ are given recursively by the Cox-de-Boor formula:
\begin{equation}
\label{order_spline}
N_{i}^k(\xi)= \displaystyle\frac{\xi - \xi_{i}}{\xi_{i+k} - \xi_{i}} N_{i}^{k-1}(\xi) + \displaystyle\frac{\xi_{i+k+1} - \xi}{\xi_{i+k+1} - \xi_{i+1}} N_{i+1}^{k-1}(\xi),
\end{equation}	 

\noindent in which fractions of the forma 0/0 are considered as zero. For more details, we refer to the reader to~\cite{PiegTill96}. 

For higher spatial dimensions, that is $d>1$, both parameter space and basis functions are built by tensorization. For instance, in the two-dimensional case, that we consider in this work, the parametric space is given by 
$$
\Xi_{p,h} \times \Xi_{p,h} = \{(\xi,\eta), \xi \in \Xi_{p,h}, \eta \in \Xi_{p,h}\}.
$$ 
\noindent Note that for the sake of simplicity we are taking the same spline degree $p$ and mesh size $h$ for both directions, but this is not restrictive. Hence, a bivariate B-spline basis function $N_{i,j}^p$ is given by means of tensor product of two univariate B-spline basis functions:
$$
N_{i,j}^p(\xi,\eta) = (N_i^p \otimes N_j^p)(\xi,\eta) = N_i^p(\xi) N_j^p(\eta).
$$

Thus, the spline space defined as the approximation space in \eqref{approxvariational} for the two-dimensional case of our model problem is given by

\begin{equation}\label{splinespace2}
\mathcal{S}_{h}^{p}(0,1)^2 = {\rm span} \{ N_{i,j}^p(\xi,\eta), i,j = 2, \ldots, p+m-1\}. 
\end{equation}
 


\subsection{Non-Uniform Rational B-splines}

In order to capture a wider set of complex geometries that use to appear in practice, we also introduce the so-called non-uniform rational B-splines (NURBS). Hence, by using NURBS as basis functions the full potential of IGA can be exploited. In order to construct them, a set of weights ${\omega_1, \ldots, \omega_{m+p}}$ is also needed. Then, the $i$-th univariate NURBS basis function of polynomial degree $p$ is given by
$$
R_i^p(\xi) = \frac{\omega_i N_i^p(\xi)}{\sum_{j=1}^{m+p} \omega_j N_j^p(\xi)}.
$$
In general, the two-dimensional NURBS basis functions cannot be constructed straightforwardly by tensorization since each weight is associated to each basis function. Hence, for $d=2$ a net of weights $\omega_{i,j}$ is considered and these basis functions are given by 
$$
R_{i,j}^p(\xi,\eta) = \frac{ \omega_{i,j} N_{i,j}^p(\xi,\eta)}{\sum_{k,l=1}^{m+p} \omega_{k,l} N_{k,l}^p(\xi,\eta) },
$$
\noindent with $i,j=1, \ldots, m+p$.

In both the B-splines and NURBS cases, the parameter space and the physical space $\Omega$ might not coincide since the parametric space is typically $(0,1)^d$. Thus, there is a spline transformation between them in which the control points are involved. For example, a two-dimensional parametrization of the domain ${\mathbf F}: (0,1)^2 \rightarrow \Omega$ is defined as follows
$$
{\mathbf F}(\xi,\eta) = \sum_{i=1}^{m+p} \sum_{j=1}^{m+p} {\mathbf B}_{i,j} R_{i,j}^p(\xi,\eta), 
$$
\noindent where $\lbrace {\mathbf B}_{i,j} \rbrace_{i,j= 1, \ldots, m+p}$ is a given control net. Hence, the approximation space $V_h^p$ of spline degree $p$ is finally determined by the composition of the B-splines/NURBS basis functions of spline degree $p$ with support on $(0,1)^d$ with the inverse transformation ${\mathbf F}^{-1}$, that is,
$$
V_h^p = {\rm span} \{ R_{i,j}^p \circ {\mathbf F}^{-1}, i,j = 1, \ldots, m+p\}.
$$




\section{Two-level method}\label{sec:two-level}

In this work, we propose an algebraic two-level method for solving isogeometric discretizations of an arbitrary polynomial degree in an efficient and robust way. This two-level method considers the target polynomial degree on the fine level whereas the order of the approximation at the coarse level is as low as possible, dictated by the  parametrization of the physical domain. In the following, in Section~\ref{sec:two-level1} we present the proposed two-level algorithm, specifying the components of the method. This algorithm act as a blackbox and it is purely algebraic. The problem on the coarse level can be exactly solved by using the preferred solver of the user, but it also can be approximated by using a suitable iterative method, for example using one iteration of a multigrid cycle as we will present in Section~\ref{sec:approx}. This, however, is only a choice of the authors but other possibilities can be equally valid. Finally, in Section~\ref{sec:improvement} we also show that a more aggressive coarse level can be used, improving the efficiency of the method.

\subsection{Two-level algorithm}\label{sec:two-level1}

As previously mentioned, in this section we explain the proposed algorithm and we introduce its main components. Recall that this two-level method solves an isogeometric discretization with a target polynomial degree on the fine grid by using a linear/quadratic discretization on the coarse level. Let us denote with $p$ and $p_{low}$ the polynomial orders of the discretization at the fine and coarse level respectively. A general two-level algorithm for solving the system $A_p u_p = b_p$, where $A_{p}$ denotes the isogeometric discretization of spline degree $p$, consists of the following: 
\begin{enumerate}
	\item Apply $\nu_1 $ steps of a suitable iterative method $S_p$ to the initial approximation $u_{p}^{0}$ on the fine level:
	$$u_{p}^{k} = u_{p}^{k-1} + S_p(b_p - A_{p} u_{p}^{k-1}), \quad k=1,\ldots, \nu_1.$$
	
	\item Compute the defect on the fine level $d_{p} = b_p - A_{p} u_{p}^{\nu_1}$ and restrict it to the coarse level by using the fine-to-coarse transfer operator
	$$d_{p_{low}} = I_p^{p_{low}} d_{p}. $$
	\item Compute the correction $e_{p_{low}}$ in the coarse level by solving the defect equation
	$$A_{p_{low}} e_{p_{low}} = d_{p_{low}},$$
	where $A_{p_{low}}$ denotes the isogeometric discretization of spline degree $p_{low}$.
	\item Prolongate and update the correction to the fine level by means of the coarse-to-fine transfer operator 
	$$u_{p}^{\nu_1} = u_{p}^{\nu_1} + I_{p_{low}}^p e_{p_{low}}.$$
	\item Apply $\nu_2 $ steps of the same iterative method $S_p$ to the current approximation:
	$$u_{p}^{\nu_1+k} = u_{p}^{\nu_1+k-1} + S_p(b_p - A_{p} u_{p}^{\nu_1+k-1}), \quad k=1,\ldots, \nu_2.$$
	
\end{enumerate}

It is very important, of course, the choice of the components of the algorithm. Hence, let us describe in the following the choice of the iterative method applied on the fine level, that we will call smoother, and the construction of the inter-grid transfer operators for the proposed two-level method.

\subsubsection{Smoother}

As relaxation procedure on the fine level, we propose the use of multiplicative Schwarz methods. These methods are a particular case of block-wise iterations which update simultaneously a set of unknowns at each time. They are based on  a splitting of the grid into blocks that gives rise to local problems.  There are many possibilities to construct these blocks. One can allow the blocks to overlap, giving rise to the class of overlapping block iterations, where smaller local problems are solved and combined via an additive or multiplicative Schwarz method. In this work, we consider multiplicative Schwarz iterations with maximum overlapping. Although this overlapping increases the computational cost of the method, it improves the convergence rates and thus a fewer number of iterations is required in order to reach the stopping criteria. A deep study of the computational cost of  these smoothers was presented in \cite{nuestropaper2018}.

More specifically, we can describe the multiplicative Schwarz iteration for solving the system $A_p u_p = b_p$ of size $n$ in the following way. Let us denote as $B^j_p$ the subset of unknowns involved in the $j-th$ block of size $n_p$, that is $B^j_p = \left\{u_{k_1},\ldots, u_{k_{n_p}}\right\}$ where $k_i$ is the global index of the $i-$th unknown in the block. In order to construct the matrix to solve associated with such a block, that is $A_p^{B^j_p}$, we consider the projection operator from the vector of unknowns $u_p$ to the vector of unknowns involved in the block. This results in a matrix $V_{B^j_p}$ of size $(n_p\times n)$, whose $i-$th row is the $k_i-$th row of the identity matrix of order $n$. Thus, matrix $A_p^{B^j_p}$ is obtained as $A_p^{B^j_p} = V_{B^j_p} A_p V_{B^j_p}^T$, and the iteration matrix of the multiplicative Schwarz method can be written as
$$\prod_{j=1}^{NB}\left( I-V_{B_p^j}^T(A_p^{B_p^j})^{-1}V_{B_p^j} A_p \right),$$
where $NB$ denotes the number of blocks obtained from the splitting of the grid, which corresponds to the number of small systems that have to be solved in a relaxation step of the multiplicative Schwarz smoother. 
In our particular case of maximum overlapping, $NB$ coincides with the number of grid-points and every block is related to a grid-point, involving that grid-point and its neighbors. 
In the two-dimensional case, that is the one that we deal with in this work, square blocks of size $\sqrt{n_p}\times \sqrt{n_p}$ around each grid point are considered. 
More concretely, we will use the nine-, twenty five- and forty nine-point multiplicative Schwarz smoothers, depending on the spline degree $p$. \\
Our study will be carried out up to $p=8$, but if one is interested in solving isogeometric discretizations with spline degree larger than $p=8$, only it is necessary to find the appropriate number of unknowns involved in the blocks to obtain an efficient two-level approach.\\

As it will be shown, by applying only one iteration of this smoother at the fine level we get a very simple and efficient solver. In order to obtain a robust solver with respect to the spline degree $p$, the size of the blocks will be chosen depending on the order of the discretization. In addition, we apply a three-colour version of the considered Schwarz-type smoothers since these counterparts provide, in general, better convergence rates, see \cite{nuestropaper2018}.



\subsubsection{Transfer operators}

Another important point of our two-level method is the construction of the restriction and prolongation operators. After computing the residual on the fine level, we restrict it to the coarse level by means of an $L^2$ projection among spline spaces. On the fine level, the solution of \eqref{approxvariational} is given by $u_h^{p} =\sum_{j=1}^{n_h^{p}} u_{j}^p \varphi_{j}^{p}$, where ${\rm dim } V_h^{p} = n_h^{p}$. Since the approximation of $u_h^{p} \in V_h^{p}$ is restricted by means of the restriction operator $I_p^{p_{low}}: V_h^{p} \rightarrow V_h^{p_{low}}$ to the space $V_h^{p_{low}}$, the resulting function $I_{p}^{p_{low}} u_h^{p} $ can be expanded as a linear combination of the spline basis functions of $V_h^{p_{low}}$. Consequently, there exists a vector of coefficients $u_{p_{low}} = \lbrace u_{j}^{p_{low}} \rbrace_{j=1}^{n_h^{p_{low}}}$ such that 
\begin{equation} 
\label{restriction2}
I_{p}^{p_{low}}u_h^p = \sum_{j=1}^{n_h^{p_{low}}} u_{j}^{p_{low}} \varphi_{j}^{p_{low}}.
\end{equation}
In order to obtain the relationship among the coefficients $u_{p}=\lbrace u_{j}^{p} \rbrace_{j=1}^{n_h^{p}}$ and $u_{p_{low}}$, we test both the approximation on the fine level and its restricted term with every basis function spanning $V_h^{p_{low}}$. Thus, one gets the following system of equations:
\begin{equation} 
\label{restrictionsystem}
\sum_{k=1}^{n_h^{p_{low}}} u_{k}^{p_{low}} \left( \varphi_{k}^{p_{low}} , \varphi_{i}^{p_{low}}\right) = \sum_{j=1}^{n_h^{p}} u_{j}^{p} \left(\varphi_{j}^{p} , \varphi_{i}^{p_{low}}\right), \quad  \forall i = 1, \ldots, n_h^{p_{low}}.
\end{equation}
This system can also be described as follows,
$$ M_{p_{low}}^{p_{low}} u_{p_{low}} = M_{p}^{p_{low}} u_{p}, $$
where 
$$ \left(M_{p_{low}}^{p_{low}}\right)_{i,j} = \int_{\Omega} \varphi_{i}^{p_{low}} \varphi_{j}^{p_{low}}  {\rm d} x, \qquad \left(M_{p}^{p_{low}}\right)_{i,j} = \int_{\Omega} \varphi_{i}^{p_{low}} \varphi_{j}^{p}  {\rm d} x .$$
Therefore, the restriction operator is given by $I_{p}^{p_{low}} = \left(M_{p_{low}}^{p_{low}}\right)^{-1} M_{p}^{p_{low}} $. Moreover, the prolongation operator is taken as its adjoint, that is, $I_{p_{low}}^{p} = \left(M_{p}^{p_{low}}\right)^T \left(M_{p_{low}}^{p_{low}}\right)^{-T}$. At this point, it is desirable to approximate $\left(M_{p_{low}}^{p_{low}}\right)^{-1}$ by row-sum lumping in order to avoid the computation of this inverse matrix exactly. 


Once introduced the components of the method, one iteration of our two-level algorithm is described in Algorithm~\ref{twolevelalgorithm}.

\begin{algorithm}[H]
\caption{ \textbf{: Two-level algorithm: ${\mathbf{u_{p}^0 \rightarrow u_{p}^{1}}}$}}\label{twolevelalgorithm}
\vspace{0.3cm}
\begin{algorithmic}
\STATE {\hspace{-0.5cm}
\begin{tabular}{lr}
\\ [-3.3ex]
$u_{p}^{1} = u_{p}^{0} + S_p(b_p - A_{p} u_{p}^{0})$ &  \hspace{-2cm} Apply \textbf{one step} \\ &  \hspace{-2cm} of the \textbf{multiplicative Schwarz method} on the fine level.\\
\\ [-3.3ex]
$d_{p} = b_p - A_{p} u_{p}^{1}$ &  \hspace{-2cm} Compute the defect on the fine level. \\
\\ [-3.3ex]
$d_{p_{low}} = I_p^{p_{low}} d_{p}$ &  \hspace{-2cm} Restrict the defect to the coarse level. \\
\\ [-3.3ex]
$A_{p_{low}} e_{p_{low}} = d_{p_{low}}$ &  \hspace{-2cm} Compute the correction $e_{p_{low}}$ in the coarse level \\
[0.5ex]
&  \hspace{-2cm} by solving the defect equation. \\
\\ [-3.3ex]
$u_{p}^{1} = u_{p}^{1} + I_{p_{low}}^p e_{p_{low}}$ & \hspace{-2cm}Prolongate and update the correction to the fine level. \\
\\ [-3.3ex]
\end{tabular}}
\end{algorithmic}
\vspace{0.3cm}
\end{algorithm}

Notice that it results in a very simple algorithm since only one single iteration of a multiplicative overlapping Schwarz method is applied on the fine level.

\subsection{Aproximation of the coarse level problem}\label{sec:approx}

Although there is an open choice for the solver at the coarse level, instead of solving exactly the coarse problem, it can also be approximated by using a suitable iterative method. In this work, we apply one $V(1,1)-$cycle that uses a red-black Gauss-Seidel iteration as smoother. 
Our numerical experiments show that one iteration of such an $h$-multigrid method is enough to ensure a good convergence rate. 
This choice will be theoretically supported by a suitable local Fourier analysis, which will be explained in Section~\ref{sec:lfa}.


\subsection{Improvement of the algorithm}\label{sec:improvement}

A further improvement of the algorithm can be achieved by using a more aggressive coarsening strategy. More concretely, we can take a discretization with $p_{low}$ and a mesh size $H=2h$ as the coarse level. Thus, the computational cost is reduced and the performance of the solver is improved without any significant effect on the convergence factors. Again, local Fourier analysis is able to theoretically support this approach, as we will see in Section~\ref{sec:lfa}.



\section{Local Fourier Analysis}\label{sec:lfa}

In this section we apply a local Fourier analysis pursuing different objectives. First, we use this analysis to theoretically support the proposed two-level algorithm and in particular the choice of the size of the block for the multiplicative Schwarz iteration depending on the spline degree $p$. In addition, in order to support the use of the $h-$multigrid as approximation on the coarse level, we apply a three-grid Fourier analysis, and as it will be shown very similar convergence rates to the case of the two-level with an exact solve on the coarse level are obtained. Finally, again a two-grid LFA is used to support the improvement of the algorithm presented in Section~\ref{sec:improvement}.

\subsection{Basics of LFA}\label{sec:basicsLFA}

Local Fourier analysis (LFA) is based on the Fourier transform theory, assuming that any grid function defined on an infinite grid ${\mathcal G}_h$ can be decomposed as a ``formal'' linear combination of complex exponential functions, $\varphi_h({\mathbf{\theta} },\mathbf{x}) = e^{\imath {\bold \theta} \mathbf{x}/h}$ with $\mathbf{x}\in{\mathcal G}_h$ and ${\bold \theta}\in\Theta:=(-\pi,\pi]^d$, known as Fourier modes. In particular such decomposition of the error function is considered and LFA  studies how the operators involved in the multilevel method act on these Fourier components, and in particular on the so-called Fourier space ${\mathcal F}({\mathcal G}_h):=\hbox{span}\{\varphi_h({\bold \theta},\mathbf{x})\,|\, {\bold \theta}\in\Theta\}$.

Here, we study the two level method previously introduced by using this analysis. With this purpose, we define the error propagation operator of the two level method, $T_p^{p_{low}} $, which relates the error in the iteration $m+1$, $e^{m+1}$, with the error in the previous iteration, $e^m$, that is, 
\begin{equation} 
e^{m+1} = T_p^{p_{low}} e^m = (I-I_{p_{low}}^p A_{p_{low}}^{-1} I^{p_{low}}_p A_p) S_pe^m.
\end{equation}
In the previous expression, $A_p$ and $A_{p_{low}}$ correspond to the IGA discrete operators of order $p$ and $p_{low}$, respectively; $I_{p_{low}}^p$ and $I^{p_{low}}_p$ are the inter-grid transfer operators, and $S_p$ represents the multiplicative Schwarz iteration which is applied within the two level method.
It is easy to see that the Fourier modes are eigenfunctions of all the operators involved in the two level method. Notice that, in this case, the transfer operators between levels do not couple Fourier modes unlike the inter-grid transfer operators within the standard h-multigrid method. 
Thus, the Fourier symbol of the error transfer operator for ${\mathbf{\theta}} \in \Theta$ is given by,
$$\widetilde{T}_p^{p_{low}}(\mathbf{\theta}) = (\widetilde{I}(\mathbf{\theta})-\widetilde{I}_{p_{low}}^p(\mathbf{\theta}) \widetilde{A}_{p_{low}}^{-1}(\mathbf{\theta}) \widetilde{I}^{p_{low}}_p(\mathbf{\theta}) \widetilde{A}_p(\mathbf{\theta}) )\widetilde{S}_p(\mathbf{\theta}).$$
The Fourier symbols of the multiplicative Schwarz smoothers considered in this work, $\widetilde{S}_p(\mathbf{\theta})$, can be found in \cite{nuestropaper2018}, and the symbols for the discrete operators and the transfer operators are easily obtained from their definitions. In this way, the asymptotic convergence factor of the two level method can be estimated by the following expression
\begin{equation}
\rho_{2g} = \sup_{{\mathbf{\theta}}\in \Theta} |\widetilde{T}_p^{p_{low}}({\mathbf{\theta}})|.
\end{equation}
 
In order to support the approximation approach by multigrid method on the coarse level given in Section \ref{sec:approx}, we have to take into account a smoothing effect at the second level and, by means of a standard coarsening on the mesh size $h$, a third level whose discretization corresponds to the spline space $\mathcal{S}^{2h}_{p_{low}}(0,1)^d$. Thus, a three-grid analysis is required and a smoother $S_{p_{low}}$ is considered. For this purpose, we introduce the error propagation matrix $M^{p_{low},2h}_{p,h}$ as follows:
$$ M^{p_{low},2h}_{p,h} = (I-I_{p_{low}}^p (I-(M_{p_{low},h}^{p_{low},2h})) A^{-1}_{p_{low}} I_{p}^{p_{low}} A_p)S_p ,$$
where $M_{p_{low},h}^{p_{low},2h}$ is the two-grid operator between the second and third levels, that is,
$$M_{p_{low},h}^{p_{low},2h} = S^{\nu_{2}}_{p_{low}} (I-I_{2h}^h  A^{-1}_{p_{low},2h} I_{h}^{2h} A_{p_{low}})S^{\nu_{1}}_{p_{low}},$$
with $I^h_{2h}$ and $I^{2h}_h$ the standard inter-grid transfer operators between the grids of size $h$ and $2h$. In addition, $\nu_{1}$ and $\nu_{2}$ denote the number of pre- and post-smoothing steps of the smoother $S_{p_{low}}$ on the second level.

\noindent In this case, in the transition from the second to the third level, some Fourier modes are coupled. Hence, we split the Fourier components into high- and low-frequency components on ${\mathcal G}_h$. The low-frequency Fourier components are those associated with frequencies belonging to $\Theta_{2h} = (-\pi/2,\pi/2]^d$. Thus, each low-frequency $\theta^0 = \theta^{00} = (\theta^{00}_1, \theta_2^{00})\in \Theta_{2h}$ is coupled with three high frequencies $\theta^{11}$, $\theta^{10}$, $\theta^{01}$, given by $\theta^{ij} = \theta^{00} - (i \, \textit{sign}(\theta^{00}_1), j\, \textit{sign}(\theta^{00}_2))\pi,\; i,j=0,1$, giving rise to the so-called spaces of $2h-$harmonics: $$\mathcal{F}^2(\theta^{00}) = \textit{span}\left\lbrace \varphi_h(\theta^{00},\cdot), \varphi_h(\theta^{11},\cdot),\varphi_h(\theta^{10},\cdot), \varphi_h(\theta^{01},\cdot)\right\rbrace, \; \hbox{with} \; \theta^{00}\in\Theta_{2h}.$$ 
Based on this decomposition of the Fourier space in terms of the subspaces of $2h-$harmonics, the spectral radius of the three-grid operator can be computed as follows:
$$\rho_{3g} = \rho (M_{p,h}^{p_{low},2h}) = \sup_{{\mathbf{\theta^{00}}}\in \Theta_{2h}} \rho (\widetilde{M}_{p,h}^{p_{low},2h}(\theta^{00})).$$

Finally, in order to analyze the improved version of the two-grid algorithm given in Section~\ref{sec:improvement}, we apply a two-grid LFA in which from the fine to the coarse levels we reduce the polynomial degree from $p$ to $p_{low}$ and also we double the grid-size from $h$ to $2h$. This two-grid analysis couples Fourier modes as explained before, and the corresponding error transfer operator is given by:
$$T_{p,h}^{p_{low},2h} = (I-I_{p_{low},2h}^{p,h} A_{p_{low},2h}^{-1} I^{p_{low},2h}_{p,h} A_{p}) S_p,$$
where the transfer operators between the spaces $\mathcal{S}_{p}^{h}(0,1)^d$ and $\mathcal{S}_{p_{low}}^{2h}(0,1)^d$, that is $I_{p_{low},2h}^{p,h}$ and $I^{p_{low},2h}_{p,h}$, are obtained by composition of $I_p^{p_{low}},I^p_{p_{low}}$ and the transfer operators $I_h^{2h}$, $I_{2h}^{h}$ between spline spaces with equal spline degree but different mesh size $h$ and $2h$. 
From this expression, the asymptotic convergence factor of the improved two level method can be estimated by the following expression:
\begin{equation}\label{2g_aggressive}
\rho_{2g}^{ag} = \sup_{{\mathbf{\theta^{00}}}\in \Theta_{2h}} \rho (\widetilde{T}_{p,h}^{p_{low},2h}(\theta^{00})).
\end{equation}

\subsection{Local Fourier analysis results}\label{sec:LFAresults}

Next, we show some LFA results to demonstrate the good performance of the proposed two level method. Firstly, we consider a linear discretization as the second level, that is, $p_{low}=1$. In Table \ref{table_LFA_Schwarz}, the two-level convergence factors predicted by LFA, $\rho_{2g}$, are shown together with the asymptotic convergence factors, $\rho_h$, obtained numerically for different values of the spline degree $p$ varying from $p=2$ to $p=8$. The asymptotic converge factors are obtained numerically by solving problem~\eqref{modelproblem} with a zero right-hand side and a random initial guess. We consider the $9-$point, $25-$point and $49-$point multiplicative Schwarz iterations at the first level. It can be seen from Table~\ref{table_LFA_Schwarz} that the factors predicted by LFA match very accurately the asymptotic convergence factors numerically obtained, and therefore the LFA results in a very useful tool to analyze the performance of the method. 
It is also observed from the table that choosing an appropriate multiplicative Schwarz smoother for each polynomial degree $p$, we obtain a robust solver with respect to $p$. This choice of the size of the blocks in the relaxation depending on the spline degree is done taking into account the two-grid convergence factors provided by the LFA, as well as the computational cost of the algorithm. In particular, we choose blocks of size $3 \times 3$ ($9-$point Schwarz smoother) for the cases $p=2, 3, 4$, blocks of size $5\times 5$  ($25-$point Schwarz smoother) for the cases $p=5, 6$ and  blocks of size $7\times 7$  ($49-$point Schwarz smoother) for spline degree $p= 7, 8$. For a more detailed explanation of how to choose the size of the blocks of the multiplicative Schwarz relaxations for different values of $p$, in terms of the LFA results and the computational cost, we refer the reader to~\cite{nuestropaper2018}.

\begin{table}[htb]
	\begin{center}
		\begin{tabular}{ccccccc}
			\cline{2-7}
			& \multicolumn{2}{c}{9p Schwarz} & \multicolumn{2}{c}{25p Schwarz} & \multicolumn{2}{c}{49p Schwarz} \\
			\cline{2-7}
			& $\rho_{2g}$ & $\rho_h$  & $\rho_{2g}$ & $\rho_h$ & $\rho_{2g}$ & $\rho_h$ \\
			\hline
			$p = 2$ &  0.1234 &  0.1212  & 0.0813 & 0.0752 & 0.0604 & 0.0725 \\
			$p = 3$ &  0.2150 &  0.2141  & 0.0874 & 0.0854 & 0.0622 & 0.0712\\
			$p = 4$ &  0.4581 &  0.4558  & 0.1294 & 0.1466 & 0.0697 & 0.0852\\
			$p = 5$ &  0.7095 &  0.7058  & 0.2690 & 0.2847 & 0.1001 & 0.1215\\
			$p = 6$ &  0.8786 &  0.8756  & 0.4549 & 0.4555 & 0.1909 & 0.2113\\
			$p = 7$ &  0.9576 &  0.9573  & 0.6623 & 0.6601 & 0.3260 & 0.3284\\
			$p = 8$ &  0.9868 &  0.9851  & 0.8278 & 0.8146 & 0.4885 & 0.4764\\
			\hline
		\end{tabular}
		\caption{Two-level ($\rho_{2g}$) convergence factors predicted by LFA together with the asymptotic convergence factors obtained numerically ($\rho_h$), for different values of the spline degree $p$. In this case, the second level is a linear discretization with the same mesh size $h$ considered for the first level.}
		\label{table_LFA_Schwarz}
	\end{center}
\end{table}

Next, we present some LFA results in order to support the approach proposed in Section~\ref{sec:approx}. In this case, one single iteration of a $V(1,1)-$cycle using red-black Gauss-Seidel as smoother is considered to approximate the problem on the coarse level. Thus, in order to analyze such approximation, we need to use the three-grid local Fourier analysis introduced in Section~\ref{sec:basicsLFA}. In Table~\ref{table_LFA_Schwarz_3g}, we show the three-grid convergence factors ($\rho_{3g}$) provided by LFA. One can observe that the predictions provided by the three-grid LFA match very well with the two-grid convergence factors predicted by the analysis for the two-level algorithm (with exact solve on the coarse level) shown in Table~\ref{table_LFA_Schwarz}. 

\begin{table}[htb]
	\begin{center}
		\begin{tabular}{cccc}
			\cline{2-4}
			& 9p Schwarz & 25p Schwarz & 49p Schwarz \\
			\hline
			$p = 2$ & 0.1281 & 0.0847 & 0.0624 \\
			$p = 3$ & 0.2144 & 0.0920 & 0.0690 \\
			$p = 4$ & 0.4566 & 0.1290 & 0.0733 \\
			$p = 5$ & 0.7078 & 0.2676 & 0.0986 \\
			$p = 6$ & 0.8773 & 0.4549 & 0.1909 \\
			$p = 7$ & 0.9569 & 0.6591 & 0.3174 \\
			$p = 8$ & 0.9864 & 0.8250 & 0.4734 \\
			\hline
		\end{tabular}
		\caption{Three-level ($\rho_{3g}$) convergence factors predicted by LFA, for different values of the spline degree $p$.}
		\label{table_LFA_Schwarz_3g}
	\end{center}
\end{table} 

Finally, we want to analyze the improvement of the algorithm presented in Section~\ref{sec:improvement}. In order to do this, we need to consider that in the second level of the algorithm we now assume a grid-size $2h$ in addition of the reduction of the spline degree to $p_{low}$. 
Again, LFA is able to support this approach by using a two-grid analysis. In Table \ref{table_2D}, the two-level convergence factors provided by this analysis (see expression in~\eqref{2g_aggressive}) are shown.

\begin{table}[htb]
	\begin{center}
		\begin{tabular}{cccc}
			\cline{2-4}
			& \multicolumn{1}{c}{9p Schwarz} & \multicolumn{1}{c}{25p Schwarz} & \multicolumn{1}{c}{49p Schwarz} \\
			\hline
			$p = 2$ &  0.1723 & 0.1137 & 0.0837   \\ 
			$p = 3$ &  0.2145 & 0.1152 & 0.0863   \\ 
			$p = 4$ &  0.4566 & 0.1290 & 0.0874   \\  
			$p = 5$ &  0.7078 & 0.2676 & 0.0986   \\  
			$p = 6$ &  0.8773 & 0.4549 & 0.1909   \\  
			$p = 7$ &  0.9569 & 0.6591 & 0.3174   \\  
			$p = 8$ &  0.9864 & 0.8250 & 0.4734   \\  
			\hline
		\end{tabular}
		\caption{Two-grid ($\rho_{2g}^{ag}$) convergence factors predicted by LFA for different values of the spline degree $p$, for the improved version of the algorithm.} 
		\label{table_2D}
	\end{center}
\end{table}

Given that this last approach is more efficient and does not deteriorate the performance of the two-level method introduced before, this will be the strategy used in the numerical experiments section.

\section{Numerical experiments}\label{sec:numerical}

In order to support the robustness and efficiency of the proposed two-level method, we have considered two different numerical experiments. In the first one, we deal with a bidimensional problem on a square domain and finally we consider another bidimensional problem whose physical domain is a quarter annulus. For the first numerical experiment we consider B-splines as basis functions and $p_{low}=1$, whereas for the second numerical experiment NURBS are used in order to exactly describe the geometry for the considered domain and therefore $p_{low} = 2$ is considered. 

As it was mentioned in Section \ref{sec:two-level}, we consider only one step of the coloured version of the multiplicative Schwarz method at the fine level. Instead of solving exactly at the coarse level, we follow the approximation strategy proposed in Section \ref{sec:approx} with the improvement introduced in Section~\ref{sec:improvement}. In all the cases the initial guess is taken as a random vector and the stopping criterion for our two-level solver is set to reduce the initial residual by a factor of $10^{-8}$. All the methods have been implemented in our in-house Fortran code, and the numerical computations have been carried out on an hp pavilion laptop 15-cs0008ns with a Core i7-8550U with 1,80 GHz and 16 GB RAM,  running Windows 10.

\subsection{Square domain}

Now, let us apply our two-level method based on overlapping multiplicative Schwarz iterations on a two-dimensional problem defined on a square domain $\Omega = (0,1)^2$. We consider the following problem:
$$
\left \{
\begin{array}{l}
- \Delta u = 2 \pi^2 \sin(\pi x) \sin(\pi y), \quad (x,y) \in \Omega, \nonumber \\
u(x,y) = 0, \quad (x,y) \; {\mbox on} \; \partial \Omega .
\end{array}
\right.
$$
%
For this numerical experiment, we consider the two-dimensional version of the spline space given in \eqref{splinespace} for different degrees ranging from $p=2$ until $p=8$. In addition, we consider a linear discretization for the coarse level and the size of the blocks is chosen depending on the spline degree. We choose blocks of size $3 \times 3$ for the cases $p=2, 3, 4$, blocks of size $5\times 5$ for the cases $p=5, 6$ and  blocks of size $7\times 7$ for spline degree $p= 7, 8$. 

\begin{table}[htbp]
	\centering
	\begin{tabular}{c|c|c|c|c|c|c|c|}  
		\cline{2-8}
		& \multicolumn{3}{|c|}{Color 9p Schwarz} & \multicolumn{2}{|c|}{Color 25p Schwarz} & \multicolumn{2}{|c|}{Color 49p Schwarz}\\
		\cline{2-8}
		& \multicolumn{1}{|c|}{$p=2$} & \multicolumn{1}{|c|}{$p=3$} & \multicolumn{1}{|c|}{$p=4$} & \multicolumn{1}{|c|}{$p=5$} & \multicolumn{1}{|c|}{$p=6$} & \multicolumn{1}{|c|}{$p=7$} & \multicolumn{1}{|c|}{$p=8$} \\
		\hline
		\multicolumn{1}{|c|}{Grid}  & $it \quad cpu$ &  $it \quad cpu$ & $it \quad cpu$ & $it \quad cpu$ & $it \quad cpu$ & $it \quad cpu$ & $it \quad cpu$ \\
		\hline
		\multicolumn{1}{|c|}{$64^2$} & $ 6 \quad 0.05$ & $ 6 \quad 0.06$ & $ 7 \quad 0.08$ & $ 5 \quad 0.19$ & $ 5 \quad 0.25$ & $ 4 \quad 0.89$ & $ 4 \quad 1.04 $ \\ 
		\multicolumn{1}{|c|}{$128^2$} & $ 6 \quad 0.14$ & $ 6 \quad 0.17$ & $ 7 \quad 0.23$ & $ 5 \quad 0.50$ & $ 5 \quad 0.62$ & $ 4 \quad 2.14$ & $ 4 \quad 2.50$ \\ 
		\multicolumn{1}{|c|}{$256^2$} & $ 6 \quad 0.46$ & $ 6 \quad 0.55$ & $ 7 \quad 0.79$ & $ 5 \quad 1.45$ & $ 5 \quad 1.78$ & $ 5 \quad 6.68$ & $ 5 \quad 7.63$ \\ 
		\multicolumn{1}{|c|}{$512^2$} & $ 6 \quad 1.69$ & $ 6 \quad 2.06$ & $ 7 \quad 2.86$ & $ 5 \quad 4.87$ & $ 5 \quad 5.83$ & $ 5 \quad 18.14 $ & $ 5 \quad 20.71$ \\ 
		\multicolumn{1}{|c|}{$1024^2$} & $ 6 \quad 6.83$ & $ 6 \quad 8.34$ & $ 7 \quad 11.21$ & $ 6 \quad 21.22$ & $ 5 \quad 21.12$ & $ 5 \quad 55.76$ & $ 5 \quad 62.86$ \\
		\hline
	\end{tabular}
	\caption{Square domain problem. Number of the proposed two-level method iterations ($it$) and computational time ($cpu$) necessary to reduce the initial residual in a factor of $10^{-8},$ for different mesh-sizes $h$ and for different values of the spline degree $p$, using the most appropriate coloured multiplicative Schwarz smoother for each $p$.}	
	\label{tab:Sch_2D}
\end{table}

In Table \ref{tab:Sch_2D}, we show the number of iterations ($it$) and the cpu time ($cpu$) in seconds needed to reach the stopping criterion for several mesh sizes and different spline degrees $p = 2, \ldots, 8$. We observe that in both cases the iteration numbers are robust with respect to the size of the grid $h$ and the spline degree $p$. With these results, we can conclude that our
two-level method provides an efficient and robust solver for B-spline isogeometric discretizations.

\subsection{Quarter annulus}

For the last experiment, our goal is to apply the two-level method to a two-dimensional problem defined in a nontrivial geometry. Thus, we set as physical domain the quarter of an annulus, 
$$
\Omega = \{(x,y) \in \mathbb{R}^2 \, | \, r^2 \leq x^2 + y^2 \leq R^2, x, y \geq 0\},
$$
where $r=0.3$, $R=0.5$. Hence, we consider the solution of the Poisson problem in such domain with homogeneous Dirichlet boundary conditions
$$
\left \{
\begin{array}{l}
- \Delta u = f(x,y), \quad (x,y) \in \Omega, \nonumber \\
u(x,y) = 0, \quad (x,y) \; {\mbox on} \; \partial \Omega, \nonumber
\end{array}
\right.
$$
where $f(x,y)$ is such that the exact solution is 
$$
u(x,y) = \sin(\pi x) \sin(\pi y) (x^2+y^2-r^2)(x^2+y^2-R^2).
$$

In order to construct this computational domain, the use of quadratic NURBS basis funcions is required. Thus, we consider discretizations of degree $p=3,\ldots,8$ with maximal smoothness for the fine level whereas the quadratic discretization is used at the coarse level. In this case, we compare the performance of the multigrid method (MG) proposed in \ref{sec:improvement} with a two-level based on a direct solver (DS) at the second level. For this purpose, in Table \ref{tab:annulus} we show the number of iterations needed to reach the stopping criterion for several mesh sizes and different spline degrees $p = 3, \ldots, 8$. We observe that the use of the mentioned MG at the coarse level slightly increases the number of iterations for some cases. Finally, we conclude that our two-level method provides an efficient and robust solver also for NURBS discretizations.

\begin{table}[htbp]
	\centering
	\begin{tabular}{c|c|c|c|c|c|c|}  
		\cline{2-7}
		& \multicolumn{2}{|c|}{C. 9p Schwarz} & \multicolumn{2}{|c|}{C. 25p Schwarz} & \multicolumn{2}{|c|}{C. 49p Schwarz}\\
		\cline{2-7} 
		& $p=3$ & $p=4$ & $p=5$ & $p=6$ & $p=7$ & $p=8$  \\
		\hline
		\multicolumn{1}{|l|}{Grid}  & DS \, MG & DS \, MG & DS \, MG & DS \, MG & DS \, MG & DS \, MG \\
		\hline
		\multicolumn{1}{|l|}{$32^2$}  & $5 \quad 5$ & $8 \quad  8 $ & $4 \quad  4 $ & $6 \quad  6 $ & $3 \quad  3  $ & $4 \quad  4  $ \\ 
		\multicolumn{1}{|l|}{$64^2$} & $5 \quad  7 $ & $8 \quad  8 $ & $4 \quad  5 $ & $6 \quad  6 $ & $4 \quad  4 $ & $5 \quad  5$ \\ 
		\multicolumn{1}{|l|}{$128^2$} & $6 \quad  8 $ & $7 \quad  8 $ & $4 \quad  6 $ & $6 \quad  6 $ & $4 \quad  5 $ & $5 \quad  5$ \\ 
		\multicolumn{1}{|l|}{$256^2$} & $6 \quad  9 $ & $ 8 \quad 8$ & $4 \quad  6$ & $6 \quad  6 $ & $4 \quad  6 $ & $ 5 \quad 6 $ \\ 
		\hline
	\end{tabular}
	\caption{Quarter annulus problem. Number of the proposed two-level method iterations ($it$) necessary to reduce the initial residual in a factor of $10^{-8},$ for different mesh-sizes $h$ and for different values of the spline degree $p$, using the most appropriate coloured multiplicative Schwarz smoother for each $p$.}	\label{tab:annulus}	
\end{table}

\section{Conclusions} \label{sec:conclusions}

In this work, we propose a purely algebraic two-level method for solving isogeometric discretizations of an arbitrary polynomial degree in an efficient and robust way. The algorithm considers the target polynomial degree on the fine level and a linear or quadratic approximation on the coarse level dictated by the parametrization of the physical domain. The method acts as a blackbox in which only one iteration of an appropriate multiplicative Schwarz method is applied on the fine level, and the coarse level can be exactly solved by using well-known techniques for solving linear and quadratic discretizations. The user can choose the preferred approach on the coarse level, but here we propose to approximate the coarse problem by using one single iteration of a suitable $h-$multigrid. In particular, we apply one $V(1,1)-$cycle based on a red-black Gauss-Seidel smoother. An enhancement of the performance of the solver is obtained if we apply a standard coarsening strategy from the first to the second level by considering a grid of size $h$ on the fine level and a coarse grid-size of $2h$. The good convergence results of the proposed method are theoretically supported by two- and three-grid local Fourier analysis and also they are demonstrated by means of two numerical experiments. 


\bibliography{mybibfile}

\begin{thebibliography}{10}
\expandafter\ifx\csname url\endcsname\relax
  \def\url#1{\texttt{#1}}\fi
\expandafter\ifx\csname urlprefix\endcsname\relax\def\urlprefix{URL }\fi
\expandafter\ifx\csname href\endcsname\relax
  \def\href#1#2{#2} \def\path#1{#1}\fi

\bibitem{Hughes_CMAME2005}
T.~Hughes, J.~Cottrell, Y.~Bazilevs, Isogeometric analysis: {CAD}, finite
  elements, {NURBS}, exact geometry and mesh refinement, Computer Methods in
  Applied Mechanics and Engineering 194~(39) (2005) 4135 -- 4195.
\newblock \href {http://dx.doi.org/https://doi.org/10.1016/j.cma.2004.10.008}
  {\path{doi:https://doi.org/10.1016/j.cma.2004.10.008}}.

\bibitem{Collier2012}
N.~Collier, D.~Pardo, L.~Dalcin, M.~Paszynski, V.~Calo, The cost of continuity:
  A study of the performance of isogeometric finite elements using direct
  solvers, Computer Methods in Applied Mechanics and Engineering 213-216 (2012)
  353 -- 361.
\newblock \href {http://dx.doi.org/https://doi.org/10.1016/j.cma.2011.11.002}
  {\path{doi:https://doi.org/10.1016/j.cma.2011.11.002}}.

\bibitem{Collier2013}
N.~O. Collier, L.~Dalc{\'i}n, D.~Pardo, V.~M. Calo, The cost of continuity:
  performance of iterative solvers on isogeometric finite elements, SIAM J.
  Scientific Computing 35.

\bibitem{Buffa2013}
A.~Buffa, H.~Harbrecht, A.~Kunoth, G.~Sangalli, {BPX}-preconditioning for
  isogeometric analysis, Computer Methods in Applied Mechanics and Engineering
  265 (2013) 63 -- 70.
\newblock \href {http://dx.doi.org/https://doi.org/10.1016/j.cma.2013.05.014}
  {\path{doi:https://doi.org/10.1016/j.cma.2013.05.014}}.

\bibitem{Beirao2012}
L.~Beir\~ao~da Veiga, D.~Cho, L.~F. Pavarino, S.~Scacchi, Overlapping {S}chwarz
  methods for isogeometric analysis, SIAM Journal on Numerical Analysis 50~(3)
  (2012) 1394--1416.
\newblock \href {http://dx.doi.org/10.1137/110833476}
  {\path{doi:10.1137/110833476}}.

\bibitem{Beirao2013}
L.~Beir\~ao~da Veiga, D.~Cho, L.~F. Pavarino, S.~Scacchi, {BDDC}
  preconditioners for isogeometric analysis, Mathematical Models and Methods in
  Applied Sciences 23~(06) (2013) 1099--1142.
\newblock \href {http://dx.doi.org/10.1142/S0218202513500048}
  {\path{doi:10.1142/S0218202513500048}}.

\bibitem{Gahalaut2013}
K.~Gahalaut, S.~Tomar, J.~Kraus, Algebraic multilevel preconditioning in
  isogeometric analysis: Construction and numerical studies, Computer Methods
  in Applied Mechanics and Engineering 266 (2013) 40 -- 56.
\newblock \href {http://dx.doi.org/https://doi.org/10.1016/j.cma.2013.07.002}
  {\path{doi:https://doi.org/10.1016/j.cma.2013.07.002}}.

\bibitem{Sangalli2016}
G.~Sangalli, M.~Tani, Isogeometric preconditioners based on fast solvers for
  the {S}ylvester equation, SIAM Journal on Scientific Computing 38~(6) (2016)
  A3644--A3671.
\newblock \href {http://dx.doi.org/10.1137/16M1062788}
  {\path{doi:10.1137/16M1062788}}.

\bibitem{Donatelli2015}
M.~Donatelli, C.~Garoni, C.~Manni, S.~Serra-Capizzano, H.~Speleers, Robust and
  optimal multi-iterative techniques for {I}g{A} {G}alerkin linear systems,
  Computer Methods in Applied Mechanics and Engineering 284 (2015) 230 -- 264,
  {I}sogeometric Analysis Special Issue.
\newblock \href {http://dx.doi.org/https://doi.org/10.1016/j.cma.2014.06.001}
  {\path{doi:https://doi.org/10.1016/j.cma.2014.06.001}}.

\bibitem{Hofreiher2015}
C.~Hofreither, W.~Zulehner, Mass smoothers in geometric multigrid for
  isogeometric analysis, in: J.-D. Boissonnat, A.~Cohen, O.~Gibaru, C.~Gout,
  T.~Lyche, M.-L. Mazure, L.~L. Schumaker (Eds.), Curves and Surfaces, Springer
  International Publishing, Cham, 2015, pp. 272--279.

\bibitem{Hofreither2017}
C.~Hofreither, S.~Takacs, W.~Zulehner, A robust multigrid method for
  isogeometric analysis in two dimensions using boundary correction, Computer
  Methods in Applied Mechanics and Engineering 316 (2017) 22 -- 42, {S}pecial
  Issue on Isogeometric Analysis: Progress and Challenges.
\newblock \href {http://dx.doi.org/https://doi.org/10.1016/j.cma.2016.04.003}
  {\path{doi:https://doi.org/10.1016/j.cma.2016.04.003}}.

\bibitem{Hofreither2017_b}
C.~Hofreither, S.~Takacs, Robust multigrid for isogeometric analysis based on
  stable splittings of spline spaces, SIAM Journal on Numerical Analysis 55~(4)
  (2017) 2004--2024.
\newblock \href {http://dx.doi.org/10.1137/16M1085425}
  {\path{doi:10.1137/16M1085425}}.

\bibitem{Tielen}
R.~Tielen, M.~Möller, D.~Göddeke, C.~Vuik, p-multigrid methods and their
  comparison to h-multigrid methods within isogeometric analysis, Computer
  Methods in Applied Mechanics and Engineering 372 (2020) 113347.
\newblock \href {http://dx.doi.org/https://doi.org/10.1016/j.cma.2020.113347}
  {\path{doi:https://doi.org/10.1016/j.cma.2020.113347}}.

\bibitem{nuestropaper2018}
A.~P\'e de~la Riva, C.~Rodrigo, F.~J. Gaspar, A robust multigrid solver for
  isogeometric analysis based on multiplicative schwarz smoothers, SIAM Journal
  on Scientific Computing 41~(5) (2019) S321--S345.

\bibitem{Bra77}
A.~Brandt, Multi-level adaptive solutions to boundary-value problems,
  Mathematics of Computation 31~(138) (1977) 333--390.

\bibitem{Bra94}
A.~Brandt, Rigorous quantitative analysis of multigrid, {I}: Constant
  coefficients two-level cycle with {L2}-norm, SIAM Journal on Numerical
  Analysis 31~(6) (1994) 1695--1730.

\bibitem{rigorousLFA}
C.~Rodrigo, F.~J. Gaspar, L.~T. Zikatanov, On the validity of the local fourier
  analysis, Journal of Computational Mathematics 37~(3) (2018) 340--348.
\newblock \href {http://dx.doi.org/https://doi.org/10.4208/jcm.1803-m2017-0294}
  {\path{doi:https://doi.org/10.4208/jcm.1803-m2017-0294}}.

\bibitem{twolevel_Ludmil}
P.~F. Antonietti, M.~Sarti, M.~Verani, L.~T. Zikatanov, A uniform additive
  {S}chwarz preconditioned for high-order discontinuous {G}alerkin
  approximations of elliptic problems, Journal of Scientific Computing 70
  (2017) 608--630.

\bibitem{deboor}
C.~d.~Boor, A Practical Guide to Splines, Springer Verlag, New York, 1978.

\bibitem{PiegTill96}
L.~Piegl, W.~Tiller, The NURBS Book, 2nd Edition, Springer-Verlag, New York,
  NY, USA, 1996.

\end{thebibliography}

\end{document}